\documentclass[12pt]{amsart}
\input{epsf}
\usepackage{amssymb,latexsym,cite}
\topskip  \numberwithin{equation}{section}
\newtheorem{theorem}{Theorem}[section]

\newtheorem{remark}[theorem]{Remark}

\begin{document}

\pagenumbering{arabic}
\pagestyle{headings}

\newcommand{\DPB}[4]{P\beta_{#1}^{(#2)}(#3,#4)}

\title{Representing polynomials by degenerate Bernoulli polynomials}
\author{Dae San Kim}
\address{Department of Mathematics, Sogang University, Seoul 121-742, Republic of Korea}
\email{dskim@sogang.ac.kr}

\author{Taekyun Kim}
\address{Department of Mathematics, Kwangwoon University, Seoul 139-701, Republic of Korea}
\email{tkkim@kw.ac.kr}


\begin{abstract}
In this paper, we consider the problem of representing any polynomial in terms of the degenerate Bernoulli polynomials and more generally of the higher-order degenerate Bernoulli polynomials. We derive explicit formulas with the help of umbral calculus and illustrate our results with some examples.
\end{abstract}
\maketitle


\section{Introduction and preliminaries}

In [11], the following identity is obtained from formulas (see Remark 3.2) expressing any polynomial in terms of Bernoulli polynomials $B_{n}(x)$ after a slight modification:
\begin{align}
& \sum_{k=1}^{n-1}\frac{1}{2k\left(2n-2k\right)}B_{2k}\left(x\right)B_{2n-2k}\left(x\right)+\frac{2}{2n-1}B_{1}\left(x\right)B_{2n-1}\left(x\right)\label{1A}\\
= & \frac{1}{n}\sum_{k=1}^{n}\frac{1}{2k}\dbinom{2n}{2k}B_{2k}B_{2n-2k}\left(x\right)+\frac{1}{n}H_{2n-1}B_{2n}\left(x\right)\nonumber 
+\frac{2}{2n-1}B_{1}\left(x\right)B_{2n-1},\nonumber
\end{align}
where $n \ge 2$ and $H_{n}=1+\frac{1}{2}+ \cdots +\frac{1}{n}$. \par
Letting $x=0$ in \eqref{1A}, we obtain a slight variant of the well-known Miki's identity:
\begin{equation}\label{2A}
\sum_{k=1}^{n-1}\frac{1}{2k\left(2n-2k\right)}B_{2k}B_{2n-2k}
= \frac{1}{n}\sum_{k=1}^{n}\frac{1}{2k}\dbinom{2n}{2k}B_{2k}B_{2n-2k}+\frac{1}{n}H_{2n-1}B_{2n}.
\end{equation}

Letting $x=\frac{1}{2}$ in \eqref{1A}, we get the Faber-Pandharipande-Zagier (FPZ) identity:
\begin{equation}\label{3A}
\sum_{k=1}^{n-1}\frac{1}{2k\left(2n-2k\right)}\overline{B}_{2k}\overline{B}_{2n-2k}
= \frac{1}{n}\sum_{k=1}^{n}\frac{1}{2k}\dbinom{2n}{2k}B_{2k}\overline{B}_{2n-2k}+\frac{1}{n}H_{2n-1}\overline{B}_{2n},
\end{equation}
where $\overline{B}_{n}=\left(2^{1-n}-1\right)B_{n}=B_{n}\left(\frac{1}{2}\right)$,
for all $n=0,1,2,\dots.$ \par
It should be emphasized that our proof of Miki's and Faber-Pandharipande-Zagier identities
follows from the simple formulas in Remark 3.2, involving only derivatives and integrals of the given polynomials.
However, the other proofs of Miki's and FPZ identities are quite involved. Indeed, for the Miki's identity, Miki [16] uses a formula for the Fermat quotient
$\frac{a^{p}-a}{p}$ modulo $p^{2}$, Shiratani-Yokoyama [20] utilizes $p$-adic analysis and Gessel [6] exploits two different expressions for Stirling numbers of the second kind $S_{2}\left(n,k\right)$. \\
In 1998, Faber and Pandharipande found that certain conjectural relations between Hodge integrals in Gromov-Witten theory require the FPZ identity. A proof of the FPZ identity was given by Zagier in the appendix of [5].
Dunne-Schubert [4] also proves the FPZ identity by making use of the asymptotic expansion
of some special polynomials coming from the quantum field theory computations. \par

Analogous formulas to Remark 3.2  can be obtained for the representations by Euler and Genocchi polynomials. Many interesting identities have been derived by using these formulas and the one in Remark 3.2 (see [7-12]). The list in the References are far from being exhaustive. The interested reader can easily find more related papers in the literature. Also, we should mention here that there are other ways of obtaining the same result as the one in \eqref{1A}. One of them is to use Fourier series expansion of the function obtained by extending by periodicity 1 of the polynomial function restricted to the interval $[0,1)$ (see [13-15]).

The aim of this paper is to derive formulas (see Theorem 3.1) expressing any polynomial in terms of the degenerate Bernoulli polynomials (see [1]) with the help of umbral calculus (see [3,18,19,21]) and to illustrate our results with some examples. This can be generalized to the higher-order degenerate Bernoulli polynomials. Indeed, we deduce formulas of representing any polynomial in terms of the higher-order degenerate Bernoulli polynomials (see Theorems 4.1 and 4.2) again by using umbral calculus. The contribution of this paper is the derivation of such formulas which, we think, have many potential applications. We have obtained formulas for representing any polynomial by the higher-order degenerate Euler polynomials and the higher-order degenerate Genocchi polynomials which will appear elsewhere.\par

The outline of this paper is as follows. In Section 1, we recall some necessary facts that are needed throughout this paper. In Section 2, we go over umbral calculus briefly. In Section 3, we derive formulas expressing any polynomial in terms of the degenerate Bernoulli polynomials. In Section 4, we derive formulas representing any polynomial in terms of the higher-order degenerate Bernoulli polynomials. In Section 5, we illustrate our results with some examples. Finally, we conclude our paper in Section 6.

The Bernoulli polynomials $B_n(x)$ are defined by
\begin{equation}\label{4A}
\frac{t}{e^t-1}e^{xt}=\sum_{n=0}^{\infty}B_{n}(x)\frac{t^n}{n!}.
\end{equation}
When $x=0$, $B_n=B_n(0)$ are called the Bernoulli numbers. We observe that $B_n(x)=\sum_{j=0}^{n}\binom{n}{j}B_{n-j}x^j,\,\frac{d}{dx}B_n(x)=nB_{n-1}(x)$.
The first few terms of $B_n$ are given by:
\begin{align*}
&B_0=1,\,B_1=-\frac{1}{2},\,B_2=\frac{1}{6},\,B_4=-\frac{1}{30},\,B_6=\frac{1}{42},\,B_8=-\frac{1}{30},\,B_{10}=\frac{5}{66},\\
&\,B_{12}=-\frac{691}{2730},\dots; B_{2k+1}=0,\,\,(k \ge 1).
\end{align*}

More generally, for any nonnegative integer $r$, the Bernoulli polynomials $B_n^{(r)}(x)$ of order $r$
are given by

\begin{equation}\label{5A}
\bigg(\frac{t}{e^t-1}\bigg)^{r}e^{xt}=\sum_{n=0}^{\infty}B_{n}^{(r)}(x)\frac{t^n}{n!}.
\end{equation}

The Euler polynomials $E_n(x)$ are defined by
\begin{equation}\label{6A}
\frac{2}{e^t+1}e^{xt}=\sum_{n=0}^{\infty}E_{n}(x)\frac{t^n}{n!}.
\end{equation}
When $x=0$, $E_n=E_n(0)$ are called the Euler numbers. We observe that $E_n(x)=\sum_{j=0}^{n}\binom{n}{j}E_{n-j}x^j,\,\frac{d}{dx}E_n(x)=nE_{n-1}(x)$.
The first few terms of $E_n$ are given by:
\begin{align*}
&E_0=1,\,E_1=-\frac{1}{2},\,E_3=\frac{1}{4},\,E_5=-\frac{1}{2},\,E_7=\frac{17}{8},\,B_9=-\frac{31}{2},\dots;\\
&E_{2k}=0,\,\,(k \ge 1).
\end{align*}

The Genocchi polynomials $G_n(x)$ are defined by
\begin{equation}\label{7A}
\frac{2t}{e^t+1}e^{xt}=\sum_{n=0}^{\infty}G_{n}(x)\frac{t^n}{n!}.
\end{equation}
When $x=0$, $G_n=G_n(0)$ are called the Genocchi numbers. We observe that $G_n(x)=\sum_{j=0}^{n}\binom{n}{j}G_{n-j}x^j$,\,$\frac{d}{dx}G_n(x)=nG_{n-1}(x)$, and $\mathrm{deg}\,G_n(x)=n-1$, for $n \ge1$.
The first few terms of $G_n$ are given by:
\begin{align*}
&G_0=0,\,G_1=1,\,G_2=-1,\,G_4=1,\,G_6=-3,\,G_8=17,\,G_{10}=-155\\
&G_{12}=2073,\dots; G_{2k+1}=0,\,\,(k \ge 1).
\end{align*}

For any nonzero real number $\lambda$, the degenerate exponentials are given by
\begin{align}
&e_{\lambda}^{x}(t)=(1+\lambda t)^{\frac{x}{\lambda}}=\sum_{n=0}^{\infty}(x)_{n,\lambda}\frac{t^n}{n!}, \label{8A}\\
&e_{\lambda}(t)=e_{\lambda}^{1}(t)=(1+\lambda t)^{\frac{1}{\lambda}}=\sum_{n=0}^{\infty}(1)_{n,\lambda}\frac{t^n}{n!}\nonumber.
\end{align}

Carlitz [1] introduced a degenerate version of the Bernoulli polynomials $B_n(x)$, called the degenerate Bernoulli polynomials and denoted by $\beta_{n,\lambda}(x)$, which is given by
\begin{equation}\label{9A}
\frac{t}{e_{\lambda}(t)-1}e_{\lambda}^{x}(t)=\sum_{n=0}^{\infty}\beta_{n,\lambda}(x)\frac{t^n}{n!}.
\end{equation}
For $x=0$, $\beta_{n,\lambda}=\beta_{n,\lambda}(0)$ are called the degenerate Bernoulli numbers. \par
More generally, for any nonnegative integer $r$, the degenerate Bernoulli polynomials $\beta_{n,\lambda}^{(r)}(x)$ of order $r$
are given by

\begin{equation}\label{10A}
\bigg(\frac{t}{e_{\lambda}(t)-1}\bigg)^{r}e_{\lambda}^{x}(t)=\sum_{n=0}^{\infty}\beta_{n,\lambda}^{(r)}(x)\frac{t^n}{n!}.
\end{equation}

We remark that $\beta_{n,\lambda}(x) \rightarrow B_{n}(x)$, and $\beta_{n,\lambda}^{(r)}(x) \rightarrow B_{n}^{(r)}(x)$, as $\lambda$ tends to $0$.

We recall some notations and facts about forward differences. Let $f$ be any complex-valued function of the real variable $x$. Then, for any real number $a$, the forward difference $\Delta_{a}$ is given by
\begin{equation}\label{11A}
\Delta_{a}f(x)=f(x+a)-f(x).
\end{equation}
If $a=1$, then we let
\begin{equation}\label{12A}
\Delta f(x)=\Delta_{1}f(x)=f(x+1)-f(x).
\end{equation}

In general, the $n$th oder forward differences are given by
\begin{equation}\label{13A}
\Delta_{a}^{n}f(x)=\sum_{i=0}^{n}\binom{n}{i} (-1)^{n-i}f(x+ia).
\end{equation}
For $a=1$, we have
\begin{equation}\label{14A}
\Delta^{n}f(x)=\sum_{i=0}^{n}\binom{n}{i} (-1)^{n-i}f(x+i).
\end{equation}
Finally, we recall that the Stirling numbers of the second kind $S_{2}(n,k)$ can be given by means of
\begin{equation}\label{15A}
\frac{1}{k!}(e^{t}-1)^{k}=\sum_{n=k}^{\infty}S_{2}(n,k)\frac{t^{n}}{n!}.
\end{equation}

\section{Review of umbral calculus}
\vspace{0.5cm}
Here we will briefly go over very basic facts about umbral calculus. For more details on this, we recommend the reader to refer to [3,18,19,21].
Let $\mathbb{C}$ be the field of complex numbers. Then $\mathcal{F}$ denotes the algebra of formal power series in $t$ over $\mathbb{C}$, given by
\begin{displaymath}
 \mathcal{F}=\bigg\{f(t)=\sum_{k=0}^{\infty}a_{k}\frac{t^{k}}{k!}~\bigg|~a_{k}\in\mathbb{C}\bigg\},
\end{displaymath}
and $\mathbb{P}=\mathbb{C}[x]$ indicates the algebra of polynomials in $x$ with coefficients in $\mathbb{C}$. \par
Let $\mathbb{P}^{*}$ be the vector space of all linear functionals on $\mathbb{P}$. If $\langle L|p(x)\rangle$ denotes the action of the linear functional $L$ on the polynomial $p(x)$, then the vector space operations on $\mathbb{P}^{*}$ are defined by
\begin{displaymath}
\langle L+M|p(x)\rangle=\langle L|p(x)\rangle+\langle M|p(x)\rangle,\quad\langle cL|p(x)\rangle=c\langle L|p(x)\rangle,
\end{displaymath}
where $c$ is a complex number. \par
For $f(t)\in\mathcal{F}$ with $\displaystyle f(t)=\sum_{k=0}^{\infty}a_{k}\frac{t^{k}}{k!}\displaystyle$, we define the linear functional on $\mathbb{P}$ by
\begin{equation}\label{1B}
\langle f(t)|x^{k}\rangle=a_{k}. 
\end{equation}
From \eqref{1B}, we note that
\begin{equation*}
 \langle t^{k}|x^{n}\rangle=n!\delta_{n,k},\quad(n,k\ge 0), 
\end{equation*}
where $\delta_{n,k}$ is the Kronecker's symbol. \par
Some remarkable linear functionals are as follows:
\begin{align}
&\langle e^{yt}|p(x) \rangle=p(y), \nonumber \\
&\langle e^{yt}-1|p(x) \rangle=p(y)-p(0), \label{2B} \\
& \bigg\langle \frac{e^{yt}-1}{t}\bigg |p(x) \bigg\rangle = \int_{0}^{y}p(u) du.\nonumber
\end{align}
Let
\begin{equation}\label{3B}
 f_{L}(t)=\sum_{k=0}^{\infty}\langle L|x^{k}\rangle\frac{t^{k}}{k!}.
\end{equation}
Then, by \eqref{1B} and \eqref{3B}, we get
\begin{displaymath}
    \langle f_{L}(t)|x^{n}\rangle=\langle L|x^{n}\rangle.
\end{displaymath}
That is, $f_{L}(t)=L$. Additionally, the map $L\longmapsto f_{L}(t)$ is a vector space isomorphism from $\mathbb{P}^{*}$ onto $\mathcal{F}$.\par  Henceforth, $\mathcal{F}$ denotes both the algebra of formal power series  in $t$ and the vector space of all linear functionals on $\mathbb{P}$. $\mathcal{F}$ is called the umbral algebra and the umbral calculus is the study of umbral algebra. 
For each nonnegative integer $k$, the differential operator $t^k$ on $\mathbb{P}$ is defined by
\begin{equation}\label{4B}
t^{k}x^n=\left\{\begin{array}{cc}
(n)_{k}x^{n-k}, & \textrm{if $k\le n$,}\\
0, & \textrm{if $k>n$.}
\end{array}\right. 
\end{equation}
Extending \eqref{4B} linearly, any power series
\begin{displaymath}
 f(t)=\sum_{k=0}^{\infty}\frac{a_{k}}{k!}t^{k}\in\mathcal{F}
\end{displaymath}
gives the differential operator on $\mathbb{P}$ defined by
\begin{equation}\label{5B}
 f(t)x^n=\sum_{k=0}^{n}\binom{n}{k}a_{k}x^{n-k},\quad(n\ge 0). 
\end{equation}
It should be observed that, for any formal power series $f(t)$ and any polynomial $p(x)$, we have
\begin{equation}\label{6B}
\langle f(t) | p(x) \rangle =\langle 1 | f(t)p(x) \rangle =f(t)p(x)|_{x=0}.
\end{equation}
Here we note that an element $f(t)$ of $\mathcal{F}$ is a formal power series, a linear functional and a differential  operator. Some notable differential operators are as follows: 
\begin{align}
&e^{yt}p(x)=p(x+y), \nonumber\\
&(e^{yt}-1)p(x)=p(x+y)-p(x), \label{7B}\\
&\frac{e^{yt}-1}{t}p(x)=\int_{x}^{x+y}p(u) du.\nonumber
\end{align}

The order $o(f(t))$ of the power series $f(t)(\ne 0)$ is the smallest integer for which $a_{k}$ does not vanish. If $o(f(t))=0$, then $f(t)$ is called an invertible series. If $o(f(t))=1$, then $f(t)$ is called a delta series. \par
For $f(t),g(t)\in\mathcal{F}$ with $o(f(t))=1$ and $o(g(t))=0$, there exists a unique sequence $s_{n}(x)$ (deg\,$s_{n}(x)=n$) of polynomials such that
\begin{equation} \label{8B}
\big\langle g(t)f(t)^{k}|s_{n}(x)\big\rangle=n!\delta_{n,k},\quad(n,k\ge 0).
\end{equation}
The sequence $s_{n}(x)$ is said to be the Sheffer sequence for $(g(t),f(t))$, which is denoted by $s_{n}(x)\sim (g(t),f(t))$. We observe from \eqref{8B} that 
\begin{equation}\label{9B}
s_{n}(x)=\frac{1}{g(t)}p_{n}(x),
\end{equation}
where $p_{n}(x)=g(t)s_{n}(x) \sim (1,f(t))$.\par
In particular, if $s_{n}(x) \sim (g(t),t)$, then $p_{n}(x)=x^n$, and hence 
\begin{equation}\label{10B}
s_{n}(x)=\frac{1}{g(t)}x^n.
\end{equation}

It is well known that $s_{n}(x)\sim (g(t),f(t))$ if and only if
\begin{equation}\label{11B}
\frac{1}{g\big(\overline{f}(t)\big)}e^{x\overline{f}(t)}=\sum_{k=0}^{\infty}\frac{s_{k}(x)}{k!}t^{k}, 
\end{equation}
for all $x\in\mathbb{C}$, where $\overline{f}(t)$ is the compositional inverse of $f(t)$ such that $\overline{f}(f(t))=f(\overline{f}(t))=t$. \par

The following equations \eqref{12B}, \eqref{13B}, and \eqref{14B} are equivalent to the fact that  $s_{n}\left(x\right)$ is Sheffer for $\left(g\left(t\right),f\left(t\right)\right)$, for some invertible $g(t)$: 
\begin{align}
f\left(t\right)s_{n}\left(x\right)&=ns_{n-1}\left(x\right),\quad\left(n\ge0\right),\label{12B}\\
s_{n}\left(x+y\right)&=\sum_{j=0}^{n}\binom{n}{j}s_{j}\left(x\right)p_{n-j}\left(y\right),\label{13B}
\end{align}
with $p_{n}\left(x\right)=g\left(t\right)s_{n}\left(x\right),$
\begin{equation}\label{14B}
s_{n}\left(x\right)=\sum_{j=0}^{n}\frac{1}{j!}\big\langle{g\left(\overline{f}\left(t\right)\right)^{-1}\overline{f}\left(t\right)^{j}}\big |{x^{n}\big\rangle}x^{j}.
\end{equation}
If $s_{n}(x)\sim(g(t),f(t))$, then the following recurrence relation holds:
\begin{equation}\label{15B}
s_{n+1}(x)=\bigg(x-\frac{g'(t)}{g(t)}\bigg)\frac{1}{f'(t)}s_{n}(x).
\end{equation}

For $s_{n}(x)\sim(g(t),f(t))$, and $r_{n}(x)\sim(h(t),l(t))$, we have 
\begin{equation}\label{16B}
s_{n}\left(x\right)=\sum_{k=0}^{n}C_{n,k}r_{k}\left(x\right),\quad\left(n\ge0\right),
\end{equation}
where 
\begin{equation}\label{17B}
C_{n,k}=\frac{1}{k!}\bigg\langle{\frac{h\left(\overline{f}\left(t\right)\right)}{g\left(\overline{f}\left(t\right)\right)}l\left(\overline{f}\left(t\right)\right)^{k}}\Big |{x^{n}}\bigg \rangle.
\end{equation}
Let $p_{n}(x),\,q_{n}(x)=\sum_{k=0}^{n}q_{n,k}x^{k}$ be sequences of polynomials. Then the umbral composition of $q_{n}(x)$ with $p_{n}(x)$ is defined to be the sequence 
\begin{equation}\label{18B}
q_{n}({\bf{p}}(x))=\sum_{k=0}^{n}q_{n,k}p_{k}(x).
\end{equation}

\section{Representation by degenerate Bernoulli polynomials}

Our interest here is to derive formulas expressing any polynomial in terms of the degenerate Bernoulli polynomials. 

From \eqref{9A} and \eqref{8A}, we first observe that 
\begin{align}
\beta_{n,\lambda}(x) \sim \big(g(t) &=\frac{e^t -1}{f(t)}=\frac{\lambda(e^t -1)}{e^{\lambda t}-1}, f(t)=\frac{1}{\lambda}(e^{\lambda t}-1)\big),\label{1C} \\
& (x)_{n,\lambda} \sim (1, f(t)=\frac{1}{\lambda}(e^{\lambda t}-1)).\label{2C}
\end{align}
From \eqref{11A}, \eqref{2B}, \eqref{7B}, \eqref{12B}, \eqref{1C} and \eqref{2C}, we note that
\begin{align}
f(t)\beta_{n,\lambda}(x)&=n\beta_{n-1}(x)=\frac{1}{\lambda}(e^{\lambda t}-1)\beta_{n}(x)=\frac{1}{\lambda}\Delta_{\lambda}\beta_{n}(x), \label{3C}\\
&f(t)(x)_{n,\lambda}=n(x)_{n-1,\lambda}.\label{4C}
\end{align}
It is immediate to see from \eqref{9A} that
\begin{align}
\Delta\beta_{n,\lambda}(x)&=\beta_{n,\lambda}(x+1)-\beta_{n,\lambda}(x)=n(x)_{n-1,\lambda}, \label{5C}\\
&\beta_{n,\lambda}(1)-\beta_{n,\lambda}=\delta_{n,1},\label{6C}
\end{align}
where $\delta_{n,1}$ is the Kronecker's delta. \par
Now, we assume that $p(x) \in \mathbb{C}[x]$ has degree $n$, and write $p(x)=\sum_{k=0}^{n}a_{k}\beta_{k,\lambda}(x)$.
Let $h(x)=p(x+1)-p(x)=\Delta p(x)$. Then, from \eqref{4C} and \eqref{5C}, we have
\begin{align}
h(x)&=\sum_{k=0}^{n}a_{k}(\beta_{k,\lambda}(x+1)-\beta_{k,\lambda}(x))\nonumber \\
&=\sum_{k=0}^{n}a_{k}k(x)_{k-1,\lambda} \label{7C}\\
&=f(t)\sum_{k=0}^{n}a_{k}(x)_{k,\lambda}\nonumber
\end{align}
For $r \ge 1$, from \eqref{7C} and \eqref{4C} we obtain
\begin{align}
(f(t))^{r-1}h(x)&=(f(t))^{r}\sum_{k=0}^{n}a_k(x)_{k,\lambda} \label{8C} \\
&=\sum_{k=r}^{n}k(k-1) \cdots (k-r+1)a_{k}(x)_{k-r,\lambda}.\nonumber
\end{align}
Letting $x=0$ in \eqref{8C}, we finally get
\begin{equation}\label{9C}
a_{r}=\frac{1}{r!}(f(t))^{r-1}h(x)|_{x=0}=\frac{1}{r!}\langle (f(t))^{r-1}|h(x) \rangle,\,\,(r \ge 1),
\end{equation}
An alternative expression of \eqref{9C} is given by\\
As $f(t)h(x)=\frac{1}{\lambda}(e^{\lambda t}-1)h(x)=\frac{1}{\lambda} \Delta_{\lambda}h(x)$,
\begin{align}
a_r&=\frac{1}{r! \lambda^{r-1}}\Delta_{\lambda}^{r-1}h(x)|_{x=0} \nonumber\\
&=\frac{1}{r! \lambda^{r-1}}(\Delta_{\lambda}^{r-1}p(1)-\Delta_{\lambda}^{r-1}p(0))\label{10C} \\
&=\frac{1}{r! \lambda^{r-1}}\Delta_{\lambda}^{r-1}\Delta p(0).\nonumber
\end{align}
From \eqref{13A}, we have another alternative expression of \eqref{9C} which is given by
\begin{align}
a_r&=\frac{1}{r! \lambda^{r-1}}\Delta_{\lambda}^{r-1}h(x)|_{x=0} \nonumber\\
&=\frac{1}{r! \lambda^{r-1}}\sum_{k=0}^{r-1}\binom{r-1}{k}(-1)^{r-1-k}h(x+k \lambda)|_{x=0} \label{11C}\\
&=\frac{1}{r! \lambda^{r-1}}\sum_{k=0}^{r-1}\binom{r-1}{k}(-1)^{r-1-k}(p(1+k \lambda)-p(k \lambda)).\nonumber
\end{align}

By using \eqref{15A}, we obtain yet another expression of \eqref{9C} that is given by
\begin{align}
a_{r}&=\frac{1}{r!}\langle (f(t))^{r-1}|h(x) \rangle \nonumber\\
&=\frac{1}{r\lambda^{r-1}}\bigg\langle \frac{1}{(r-1)!} (e^{\lambda t}-1)^{r-1}\bigg|h(x) \bigg\rangle \label{11-1C}\\
&=\frac{1}{r\lambda^{r-1}}\bigg\langle \sum_{l=r-1}^{\infty}S_{2}(l,r-1)\frac{\lambda^{l}t^{l}}{l!}\bigg|h(x) \bigg\rangle \nonumber \\
&=\frac{1}{r}\sum_{l=r-1}^{n}S_{2}(l,r-1)\frac{\lambda^{l-r+1}}{l!}(p^{(l)}(1)-p^{(l)}(0)), \nonumber
\end{align}
where $p^{(l)}(x)=(\frac{d}{dx})^{l}p(x)$.

Now, it remains to determine $a_{0}$. We first note from \eqref{2B}, \eqref{6B}, \eqref{3C} and \eqref{6C} that 
\begin{align}
g(t)p(x)|_{x=0}&=\langle g(t) | p(x) \rangle \nonumber \\
&=\bigg\langle \frac{e^t -1}{f(t)} \bigg | p(x) \bigg\rangle \nonumber\\
&=\sum_{k=0}^{n}a_{k}\bigg\langle \frac{e^t -1}{f(t)} \bigg | \beta_{k,\lambda}(x) \bigg\rangle \label{12C}
\end{align}
\begin{align}
&=\sum_{k=0}^{n}a_{k}\bigg\langle \frac{e^t -1}{f(t)} \bigg | f(t) \frac{1}{k+1}\beta_{k+1,\lambda}(x) \bigg\rangle \nonumber\\
&=\sum_{k=0}^{n}\frac{a_{k}}{k+1}\bigg\langle e^t -1 \bigg |\beta_{k+1,\lambda}(x) \bigg\rangle\nonumber \\
&=\sum_{k=0}^{n}\frac{a_{k}}{k+1}\delta_{k+1,1}=a_{0}. \nonumber
\end{align}
We want to find more explicit expression for \eqref{12C}. As $\frac{\lambda t}{e^{\lambda t}-1}e^{xt}=\sum_{n=0}^{\infty}\lambda^{n}B_n(\frac{x}{\lambda})\frac{t^n}{n!}$, we note from \eqref{10B} that $\lambda^{n}B_n(\frac{x}{\lambda})=\frac{\lambda t}{e^{\lambda t}-1}x^n$.
To proceed further, we let $p(x)=\sum_{i=0}^{n}b_{i}x^{i}$.\par
From \eqref{2B}, \eqref{1C} and \eqref{12C}, noting that $g(t)=\frac{e^t -1}{t} \frac{\lambda t}{e^{\lambda t}-1}$, we have
\begin{align}
a_{0}&=\bigg\langle \frac{e^t -1}{t}\bigg | \frac{\lambda t}{e^{\lambda t}-1}p(x) \bigg\rangle \nonumber\\
&=\sum_{i=0}^{n}b_{i}\bigg\langle \frac{e^t -1}{t}\bigg | \lambda^{i}B_{i}\big(\frac{x}{\lambda}\big)\bigg\rangle \nonumber \\
&=\bigg\langle \frac{e^t -1}{t}\bigg |\sum_{i=0}^{n}b_{i}\lambda^{i}B_{i}\big(\frac{x}{\lambda}\big)\bigg\rangle \label{13C}\\
&=\bigg\langle \frac{e^t -1}{t}\bigg |p\big(\lambda{\bf{B}}\big(\frac{x}{\lambda}\big)\big)\bigg\rangle \nonumber\\
&=\int_{0}^{1}p\big(\lambda{\bf{B}}\big(\frac{u}{\lambda}\big)\big)\, du, \nonumber
\end{align}
where $p\big(\lambda{\bf{B}}\big(\frac{x}{\lambda}\big)\big)$ denotes the umbral composition of $p(x)$ with $\lambda^{i}B_{i}\big(\frac{x}{\lambda}\big)$, that is, it is given by 
$p\big(\lambda{\bf{B}}\big(\frac{x}{\lambda}\big)\big)=\sum_{i=0}^{n}b_{i}\lambda^{i}B_{i}\big(\frac{x}{\lambda}\big)$.
Finally, from \eqref{9C}--\eqref{11C} and \eqref{13C} we get the following theorem.

\begin{theorem}
Let $p(x) \in \mathbb{C}[x], \mathrm{deg}\, p(x)=n$. Then we have 
$p(x)=\sum_{k=0}^{n}a_k \beta_{k,\lambda}(x)$, \\
where
\begin{align*}
&a_0=g(t)p(x)|_{x=0}=\bigg\langle \frac{e^{t}-1}{t}\bigg | \frac{\lambda t}{e^{\lambda t}-1}p(x) \bigg \rangle=\int_{0}^{1}p\big(\lambda{\bf{B}}\big(\frac{u}{\lambda}\big)\big) du, \\ 
&a_{k}=\frac{1}{k!}(f(t))^{k-1}(p(x+1)-p(x))|_{x=0}\\
&=\frac{1}{k! \lambda^{k-1}}\big \langle \big(e^{\lambda t}-1\big)^{k-1} \big | p(x+1)-p(x)\big \rangle \\
&=\frac{1}{k! \lambda^{k-1}}\Delta_{\lambda}^{k-1}\Delta p(0) \\
&=\frac{1}{k! \lambda^{k-1}}\sum_{j=0}^{k-1}\binom{k-1}{j}(-1)^{k-1-j}(p(1+j \lambda)-p(j \lambda)) \\
&=\frac{1}{k}\sum_{l=k-1}^{n}S_{2}(l,k-1)\frac{\lambda^{l-k+1}}{l!}(p^{(l)}(1)-p^{(l)}(0)),
\,\,\mathrm{for}\,\, k=1,2,\dots,n,
\end{align*}
\end{theorem}
where $g(t)=\frac{\lambda(e^t -1)}{e^{\lambda t}-1}, \,\,
f(t)=\frac{1}{\lambda}(e^{\lambda t}-1)$,
and $p\big(\lambda{\bf{B}}\big(\frac{x}{\lambda}\big)\big)$ denotes the umbral composition of $p(x)$ with $\lambda^{i}B_{i}\big(\frac{x}{\lambda}\big)$.

\begin{remark}
Let $p(x) \in \mathbb{C}[x], with \,\,\mathrm{deg}\, p(x)=n$. Write $p(x)=\sum_{k=0}^{n}a_kB_k(x)$. As $\lambda$ tends to $0$, $g(t) \rightarrow \frac{e^t -1}{t},\,\, f(t) \rightarrow t,\,\,p\big(\lambda{\bf{B}}\big(\frac{x}{\lambda}\big)\big) \rightarrow p(x)$. Thus we recover from Theorem 3.1 the result obtained in \textnormal{[11]}. 
Namely, we have
\begin{align*}
a_0=\int_{0}^{1}p(t) dt,\quad a_{k}=\frac{1}{k!}(p^{(k-1)}(1)-p^{(k-1)}(0)), \,\,\mathrm{for}\,\, k=1,2,\dots,n.
\end{align*}
\end{remark}

\begin{remark} We have two other expressions for $a_{0}$. \par
(a) With $g(t)=\frac{\lambda(e^t -1)}{e^{\lambda t}-1}, f(t)=\frac{1}{\lambda}(e^{\lambda t}-1)$, from \eqref{15B} we have
\begin{equation}\label{14C}
\beta_{n+1,\lambda}(x)=\bigg(x-\frac{g'(t)}{g(t)}\bigg)\frac{1}{f'(t)}\beta_{n,\lambda}(x).
\end{equation}
By making use of \eqref{6C} and \eqref{14C}, we have
\begin{align*}
\bigg\langle e^{t}-1 \bigg|\bigg(x-\frac{g'(t)}{g(t)}\bigg)\frac{1}{f'(t)}p(x) \bigg\rangle 
&=\sum_{k=0}^{n} a_{k}\langle e^{t}-1 | \beta_{k+1,\lambda}(x) \rangle \\
&=\sum_{k=0}^{n} a_{k}(\beta_{k+1,\lambda}(1)-\beta_{k+1,\lambda})\\
&=a_{0}.
\end{align*}
(b) By substituting $x=0$ into $p(x)=\sum_{k=0}^{n}a_k \beta_{k,\lambda}(x)$, we get
\begin{equation*}
a_{0}=p(0)-a_{1}\beta_{1,\lambda}- \cdots -a_{n}\beta_{n,\lambda}.
\end{equation*}
\end{remark}

\begin{remark}
Theorems 3.1 has been applied to many polynomials in oder to obtain interesting identities for certain special polynomials and numbers. Some of the polynomials that have been considered are as follows:

(a) \begin{equation*}
\sum B_{i_1}(x) \cdots B_{i_r}(x)E_{j_1}(x) \cdots E_{j_s}(x)G_{k_1+1}(x) \cdots G_{k_t+1}(x)x^{l},
\end{equation*}
where the sum is over all nonnegative integers $ i_1,\cdots, i_r, j_1,\cdots,j_s,k_1,\cdots,k_t,l$ such that
$i_1+ \cdots+ i_r +j_1+\cdots+j_s+k_1+ \cdots+k_t+l=n$, and $r, s, t, l$ are nonnegative integers with $r+s+t \ge 1$. \\
(b) \begin{equation*}
\sum \frac{B_{i_1}(x) \cdots B_{i_r}(x)E_{j_1}(x) \cdots E_{j_s}(x)G_{k_1+1}(x) \cdots G_{k_t+1}(x)x^{l}}{i_1! \cdots i_r! j_1! \cdots j_s! (k_1+1)! \cdots(k_t+1)! l!}, \\
\end{equation*}
where the sum is over all nonnegative integers $ i_1,\cdots, i_r, j_1,\cdots,j_s,k_1,\cdots,k_t,l $ such that
$i_1+ \cdots+ i_r +j_1+\cdots+j_s+k_1+ \cdots+k_t+l=n$, and $r, s, t, l$ are nonnegative integers with $r+s+t \ge 1$. \\
(c) \begin{equation*}
\sum \frac{B_{i_1}(x) \cdots B_{i_r}(x)E_{j_1}(x) \cdots E_{j_s}(x)G_{k_1+1}(x) \cdots G_{k_t+1}(x)x^{l}}{i_1 \cdots i_r j_1 \cdots j_s (k_1+1) \cdots(k_t+1) l}, \\
\end{equation*}
where the sum is over all positive integers $ i_1,\cdots, i_r, j_1,\cdots,j_s,l$ and nonnegative integers $k_1,\cdots,k_t$ such that $i_1+ \cdots+ i_r +j_1+\cdots+j_s+k_1+ \cdots+k_t+l=n$, and $r, s, t, l$ are nonnegative integers with $r+s+t \ge 1$. 
\end{remark}

\section{Representation by higher-order degenerate Bernoulli polynomials}

Our interest here is to derive formulas expressing any polynomial in terms of the higher-order degenerate Bernoulli polynomials. \par
With $g(t)=\frac{e^t -1}{f(t)}=\frac{\lambda(e^t -1)}{e^{\lambda t}-1}, f(t)=\frac{1}{\lambda}(e^{\lambda t}-1)$, from \eqref{10A} we note that
\begin{equation*}
\beta_{n,\lambda}^{(r)}(x) \sim (g(t)^{r}, f(t)).
\end{equation*}
Also, from \eqref{12B}, we have
\begin{equation}\label{1D}
f(t)\beta_{n,\lambda}^{(r)}(x)=n\beta_{n-1,\lambda}^{(r)}(x),
\end{equation}
and from the generating function of the higher-order degenerate Bernoulli polynomials, it is immediate to see that
\begin{equation}\label{2D}
\Delta\beta_{n,\lambda}^{(r)}(x)=\beta_{n,\lambda}^{(r)}(x+1)-\beta_{n,\lambda}^{(r)}(x)=n\beta_{n-1,\lambda}^{(r-1)}(x).
\end{equation}
Now, we assume that $p(x) \in \mathbb{C}[x]$ has degree $n$, and write $p(x)=\sum_{k=0}^{n}a_{k}\beta_{k,\lambda}^{(r)}(x)$.
It is important to observe from \eqref{1D} and \eqref{2D} that
\begin{align}
g(t)\beta_{n,\lambda}^{(r)}(x)&=\frac{e^t -1}{f(t)}\beta_{n,\lambda}^{(r)}(x) \nonumber \\
&=\frac{e^t -1}{f(t)}f(t)\frac{\beta_{n+1,\lambda}^{(r)}(x)}{n+1}\label{3D}\\
&=\frac{1}{n+1}(\beta_{n+1,\lambda}^{(r)}(x+1)-\beta_{n+1,\lambda}^{(r)}(x))\nonumber \\
&=\beta_{n,\lambda}^{(r-1)}(x).\nonumber
\end{align}
Thus, from \eqref{3D} we have $g(t)^{r}\beta_{n,\lambda}^{(r)}(x)=\beta_{n,\lambda}^{(0)}(x)=(x)_{n,\lambda}$, and hence 
\begin{equation}\label{4D}
g(t)^{r}p(x)=\sum_{l=0}^{n}a_{l}\,g(t)^{r}\beta_{l,\lambda}^{(r)}(x)=\sum_{l=0}^{n}a_{l}(x)_{l,\lambda}.
\end{equation}
By using \eqref{4D} and \eqref{4C}, we observe that
\begin{align}
f(t)^{k}g(t)^{r}p(x)&=\sum_{l=0}^{n}a_{l}f(t)^{k}(x)_{l,\lambda}\label{5D}\\
&=\sum_{l=k}^{n}a_{l}\,l(l-1) \cdots (l-k+1)(x)_{l-k,\lambda}.\nonumber
\end{align}
By evaluating \eqref{5D} at $x=0$, we obtain
\begin{equation}\label{6D}
a_k=\frac{1}{k!}f(t)^{k}g(t)^{r}p(x)|_{x=0}=\frac{1}{k!}\langle f(t)^{k}g(t)^{r} | p(x) \rangle.
\end{equation}
This also follows from the observation $\langle g(t)^{r}f(t)^{k} | \beta_{l,\lambda}^{(r)}(x) \rangle=l!\,\delta_{l,k}.$ \par
To proceed further, we note that $f(t)g(t)=f(t)\frac{e^{t}-1}{f(t)}=e^{t}-1$. \par
Assume first that $r > n$. Then $r > k$, for all $k=0,1, \dots ,n$.
\begin{align}
a_k&=\frac{1}{k!}\langle f(t)^{k}g(t)^{r} | p(x) \rangle \nonumber\\
&=\frac{1}{k!}\langle (e^{t}-1)^{k}g(t)^{r-k} | p(x) \rangle \nonumber\\
&=\frac{1}{k!}\sum_{j=0}^{k}(-1)^{k-j}\binom{k}{j}\langle e^{jt} |g(t)^{r-k} p(x) \label{7D}\rangle \\
&=\frac{1}{k!}\sum_{j=0}^{k}(-1)^{k-j}\binom{k}{j}g(t)^{r-k} p(x)|_{x=j}\nonumber \\
&=\frac{1}{k!}\sum_{j=0}^{k}(-1)^{k-j}\binom{k}{j}g(t)^{r-k} p(j).\nonumber
\end{align}
Next, we assume that $r \le n$. If further $0 \le k <r$, then $a_{k}$ is the same as the expression in \eqref{7D}.
Let $ r \le k \le n$. Then we have
\begin{align}
a_k&=\frac{1}{k!}\langle f(t)^{k}g(t)^{r} | p(x) \rangle \nonumber\\
&=\frac{1}{k!}\langle f(t)^{k-r}(e^{t}-1)^{r} | p(x) \rangle \label{8D}\\
&=\frac{1}{k!}\sum_{j=0}^{r}(-1)^{r-j}\binom{r}{j}\langle e^{jt} | f(t)^{k-r} p(x) \rangle\nonumber \\
&=\frac{1}{k!}\sum_{j=0}^{r}(-1)^{r-j}\binom{r}{j} f(t)^{k-r} p(j).\nonumber
\end{align}

Summarizing the results so far, from \eqref{7D} and \eqref{8D} we obtain the following theorem.
\begin{theorem}
Let $p(x) \in \mathbb{C}[x], \mathrm{deg}\, p(x)=n$. Let $g(t)=\frac{\lambda(e^t -1)}{e^{\lambda t}-1}, f(t)=\frac{1}{\lambda}(e^{\lambda t}-1)$. Then we have \\
(a) For $ r> n$, we have
\begin{equation*}
p(x)=\sum_{k=0}^{n}\frac{1}{k!}\sum_{j=0}^{k}(-1)^{k-j}\binom{k}{j}g(t)^{r-k}p(j) \beta_{k,\lambda}^{(r)}(x).
\end{equation*}
(b) For $r \le n$, we have
\begin{align*}
p(x)&=\sum_{k=0}^{r-1}\frac{1}{k!}\sum_{j=0}^{k}(-1)^{k-j}\binom{k}{j}g(t)^{r-k}p(j) \beta_{k,\lambda}^{(r)}(x) \\
&+\sum_{k=r}^{n}\frac{1}{k!}\sum_{j=0}^{r}(-1)^{r-j}\binom{r}{j}f(t)^{k-r}p(j) \beta_{k,\lambda}^{(r)}(x). 
\end{align*}
\end{theorem}

We would like to find more explicit expressions for the results in Theorem 4.1.

First, we note from \eqref{13A} and \eqref{15A} that
\begin{align}
f(t)^{k-r}p(j)&=f(t)^{k-r}p(x)|_{x=j} \nonumber\\
&=\frac{1}{\lambda^{k-r}}\Delta_{\lambda}^{k-r}p(x)|_{x=j}\label{9D}\\
&=\frac{1}{\lambda^{k-r}}\sum_{l=0}^{k-r}\binom{k-r}{l}(-1)^{k-r-l}p(j+l\lambda).\nonumber \\
&=(k-r)!\sum_{l=k-r}^{n}S_{2}(l,k-r)\frac{\lambda^{l-k+r}}{l!}p^{(l)}(j). \nonumber
\end{align}
Let $I$ be the linear integral operator defined on $\mathbb{P}$, which is given by 
\begin{equation}\label{10D}
Iq(x)=\frac{e^{t}-1}{t}q(x)=\int_{x}^{x+1}q(u)\,du.
\end{equation}
From $\big(\frac{\lambda t}{e^{\lambda t}-1}\big)^{a}e^{xt}=\sum_{n=0}^{\infty}\lambda^{n}B_{n}^{(a)}(\frac{x}{\lambda})\frac{t^n}{n!}$, we get 
\begin{equation}\label{11D}
\lambda^{n}B_{n}^{(a)}\big(\frac{x}{\lambda}\big)=\bigg(\frac{\lambda t}{e^{\lambda t}-1}\bigg)^{a}x^{n},
\end{equation}
and hence, from \eqref{11D}, we see that
\begin{equation}\label{12D}
\bigg(\frac{\lambda t}{e^{\lambda t}-1}\bigg)^{r-k}p(x)=p\big(\lambda {\bf{B}}^{(r-k)}\big(\frac{x}{\lambda}\big)\big),
\end{equation}
where $p\big(\lambda{\bf{B}}^{(r-k)}\big(\frac{x}{\lambda}\big)\big)$ denotes the umbral composition of $p(x)$ with $\lambda^{i}B_{i}^{(r-k)}\big(\frac{x}{\lambda}\big)$.\par
Now, noting that $g(t)=\frac{e^{t}-1}{t}\frac{\lambda t}{e^{\lambda t}-1}$, from \eqref{10D} and \eqref{12D}, we have
\begin{align}
g(t)^{r-k}p(j)&=g(t)^{r-k}p(x) |_{x=j}\nonumber\\
&=\bigg(\frac{e^{t}-1}{t}\bigg)^{r-k}\bigg(\frac{\lambda t}{e^{\lambda t}-1}\bigg)^{r-k}p(x) \bigg |_{x=j}\nonumber\\
&=\bigg(\frac{e^{t}-1}{t}\bigg)^{r-k} p\big(\lambda {\bf{B}}^{(r-k)}\big(\frac{x}{\lambda}\big)\big)  \bigg |_{x=j}\label{13D}\\
&=I^{r-k} p\big(\lambda {\bf{B}}^{(r-k)}\big(\frac{x}{\lambda}\big)\big)  \bigg |_{x=j} \nonumber\\
&=\bigg \langle e^{jt} \bigg(\frac{e^{t}-1}{t}\bigg)^{r-k} \bigg| \bigg(\frac{\lambda t}{e^{\lambda t}-1} \bigg)^{r-k} p(x) \bigg \rangle.\nonumber
\end{align}
From \eqref{15A}, we see that $\big(\frac{e^{t}-1}{t}\big)^{m}=\sum_{l=0}^{\infty}S_{2}(l+m,m)\frac{m!}{(l+m)!}t^{l}$, and hence from \eqref{12D} and \eqref{13D} we get another expression for $g(t)^{r-k}p(j)$ in the following:
\begin{align*}
g(t)^{r-k}p(j)&=\bigg(\frac{\lambda t}{e^{\lambda t}-1}\bigg)^{r-k}\bigg(\frac{e^{t}-1}{t}\bigg)^{r-k}p(x) \bigg |_{x=j}\\
&=\sum_{l=0}^{n}S_{2}(l+r-k,r-k)\frac{(r-k)!}{(l+r-k)!}p^{(l)}\big(\lambda{\bf{B}}^{(r-k)}\big(\frac{x}{\lambda}\big)\big) \big|_{x=j}.
\end{align*}
Now, from Theorem 4.1, \eqref{9D} and \eqref{12D}, we finally arrive at the following theorem.
\begin{theorem}
Let $p(x) \in \mathbb{C}[x], \mathrm{deg}\, p(x)=n$. Then we have \\
(a) For $ r> n$, we have
\begin{align*}
p(x)&=\sum_{k=0}^{n}\frac{1}{k!}\sum_{j=0}^{k}(-1)^{k-j}\binom{k}{j}I^{r-k} p\big(\lambda {\bf{B}}^{(r-k)}\big(\frac{x}{\lambda}\big)\big)  \bigg |_{x=j} \beta_{k,\lambda}^{(r)}(x) \\
&=\sum_{k=0}^{n}\frac{1}{k!}\sum_{j=0}^{k}\sum_{l=0}^{n}(-1)^{k-j}\binom{k}{j}\frac{(r-k)!}{(l+r-k)!}S_{2}(l+r-k,r-k) \\
&\quad\quad\quad \times p^{(l)}\big(\lambda{\bf{B}}^{(r-k)}\big(\frac{x}{\lambda}\big)\big) \big|_{x=j} \beta_{k,\lambda}^{(r)}(x).
\end{align*}
(b) For $r \le n$, we have
\begin{align*}
p(x)&=\sum_{k=0}^{r-1}\frac{1}{k!}\sum_{j=0}^{k}(-1)^{k-j}\binom{k}{j}I^{r-k} p\big(\lambda {\bf{B}}^{(r-k)}\big(\frac{x}{\lambda}\big)\big)  \bigg |_{x=j} \beta_{k,\lambda}^{(r)}(x) \\
&+\sum_{k=r}^{n}\frac{1}{k!\lambda^{k-r}}\sum_{j=0}^{r}\sum_{l=0}^{k-r}(-1)^{k-j-l}\binom{r}{j}\binom{k-r}{l}p(j+l\lambda)\beta_{k,\lambda}^{(r)}(x) \\
&=\sum_{k=0}^{r-1}\frac{1}{k!}\sum_{j=0}^{k}\sum_{l=0}^{n}(-1)^{k-j}\binom{k}{j}\frac{(r-k)!}{(l+r-k)!}S_{2}(l+r-k,r-k) \\
&\quad\quad\quad \times p^{(l)}\big(\lambda{\bf{B}}^{(r-k)}\big(\frac{x}{\lambda}\big)\big) \big|_{x=j} \beta_{k,\lambda}^{(r)}(x) \\
&+\sum_{k=r}^{n}\frac{1}{k!}\sum_{j=0}^{r}\sum_{l=k-r}^{n}(-1)^{r-j}\binom{r}{j}(k-r)!\frac{\lambda^{l-k+r}}{l!}S_{2}(l,k-r)p^{(l)}(j)\beta_{k,\lambda}^{(r)}(x) 
\end{align*}
Here $I$ is the linear operator given by $Iq(x)=\int_{x}^{x+1}q(u)\,du,$ and $p\big(\lambda{\bf{B}}^{(r-k)}\big(\frac{x}{\lambda}\big)\big)$ denotes the umbral composition of $p(x)$ with $\lambda^{i}B_{i}^{(r-k)}\big(\frac{x}{\lambda}\big)$.
\end{theorem}

\begin{remark}
Let $p(x) \in \mathbb{C}[x], with \,\,\mathrm{deg}\, p(x)=n$. Write $p(x)=\sum_{k=0}^{n}a_kB_{k}^{(r)}(x)$. As $\lambda$ tends to $0$, $g(t) \rightarrow \frac{e^t -1}{t},\,\,f(t) \rightarrow t,\,\, p\big(\lambda {\bf{B}}^{(r-k)}\big(\frac{x}{\lambda}\big)\big) \rightarrow p(x)$. Thus, from Theorem 4.1, we recover the following result obtained in \textnormal{[8]}: \par
(a) For $r>n$, we have 
\begin{equation*} 
p(x)=\sum_{k=0}^{n}\bigg(\sum_{j=0}^{k}\frac{1}{k!}(-1)^{k-j}\binom{k}{j}I^{r-k}p(j)\bigg)B_{k}^{(r)}(x).
\end{equation*}
(b) For $r \leq n$, we have 
\begin{align*}
p(x)=&\sum_{k=0}^{r-1}\bigg(\sum_{j=0}^{k}\frac{1}{k!}(-1)^{k-j}\binom{k}{j}I^{r-k}p(j)\bigg)B_{k}^{(r)}(x)\\
&+\sum_{k=r}^{n}\bigg(\sum_{j=0}^{r}\frac{1}{k!}(-1)^{r-j}p^{(k-r)}(j)\bigg)B_{k}^{(r)}(x).
\end{align*}
Here $I$ is the integral operator $Iq(x)=\int_{x}^{x+1}q(u)\,du$. 
\end{remark}

\section{Examples}
(a) Here we illustrate Theorem 3.1, with $p(x)=B_{n}(x)$. Let $B_{n}(x)=\sum_{k=0}^{n}a_{k}\beta_{k,\lambda}(x)$.\par
Then, as $B_n(x)=\sum_{j=0}^{n}\binom{n}{j}B_{n-j}x^{j}$ and $\frac{d}{dx}\frac{1}{j+1}B_{j+1}(x)=B_{j}(x)$, we have
\begin{align}
a_{0}&=\int_{0}^{1}B_{n}\big(\lambda{\bf{B}}\big(\frac{u}{\lambda}\big)\big) du \nonumber\\
&=\sum_{j=0}^{n}\binom{n}{j}B_{n-j}\lambda^{j}\int_{0}^{1}B_{j}\big(\frac{u}{\lambda}\big) du \label{1E}\\
&=\sum_{j=0}^{n}\binom{n}{j}B_{n-j}\frac{\lambda^{j+1}}{j+1}\big(B_{j+1}\big(\frac{1}{\lambda}\big)-B_{j+1}\big), \nonumber
\end{align}
which is, recalling that $\lambda^{n}B_n(\frac{x}{\lambda})=\frac{\lambda t}{e^{\lambda t}-1}x^n$ and from Theorem 3.1, also equal to
\begin{align}
a_{0}&=\bigg\langle \frac{e^{t}-1}{t}\bigg | \frac{\lambda t}{e^{\lambda t}-1}B_{n}(x) \bigg \rangle \nonumber \\
&=\bigg\langle \frac{e^{t}-1}{t}\bigg | \frac{\lambda t}{e^{\lambda t}-1}\frac{t}{e^{t}-1}x^n \bigg  \rangle \nonumber\\
&=\bigg\langle 1 \bigg | \frac{\lambda t}{e^{\lambda t}-1}x^n \bigg  \rangle\label{2E} \\
&=\bigg\langle 1 \bigg | \lambda^{n}B_{n}\big(\frac{x}{\lambda}\big) \bigg  \rangle \nonumber\\
&=\lambda^{n}B_{n}. \nonumber
\end{align}
Incidentally, \eqref{1E} and \eqref{2E} together give us the following polynomial identity:
\begin{equation*}
\sum_{j=0}^{n}\binom{n}{j}B_{n-j}\frac{x^{j+1}}{j+1}\big(B_{j+1}\big(\frac{1}{x}\big)-B_{j+1}\big)=x^{n}B_{n}.
\end{equation*}
Further, for $k$ with $1 \le k \le n$, we have
\begin{align}
a_{k}&=\frac{1}{k! \lambda^{k-1}}\big \langle \big(e^{\lambda t}-1\big)^{k-1} \big | B_n(x+1)-B_n(x)\big \rangle \nonumber\\
&=\frac{1}{k! \lambda^{k-1}}\big \langle \big(e^{\lambda t}-1\big)^{k-1} \big | nx^{n-1}\big \rangle \label{3E}\\
&=\frac{1}{k! \lambda^{k-1}}\Delta_{\lambda}^{k-1}nx^{n-1}|_{x=0}\nonumber\\
&=\frac{n}{k! \lambda^{k-1}}\Delta_{\lambda}^{k-1}0^{n-1}.\nonumber
\end{align}
Noting that $B_n(x) \sim \big(\frac{e^{t}-1}{t},t \big),\,\, \beta_{n,\lambda}(x) \sim \big(\frac{\lambda(e^t -1)}{e^{\lambda t}-1}, \frac{1}{\lambda}(e^{\lambda t}-1)\big)$, one can determine the coefficients also by using \eqref{17B}. Thus, from \eqref{2E} and \eqref{3E} we have 
\begin{equation*}
B_{n}(x)=\lambda^{n}B_{n}\beta_{0,\lambda}(x)+\sum_{k=1}^{n}\frac{n}{k! \lambda^{k-1}}\Delta_{\lambda}^{k-1}0^{n-1}\beta_{k,\lambda}(x).
\end{equation*}
(b) Here we illustrate Theorem 3.1, for $p(x)=\sum_{k=1}^{n-1}\frac{1}{k(n-k)}B_{k}(x)B_{n-k}(x),\,\,(n \ge 2)$. For this, we first recall from [11] that
\begin{equation}\label{4E}
p(x)=\frac{2}{n}\sum_{l=0}^{n-2}\frac{1}{n-l}\binom{n}{l}B_{n-l}B_l(x)+\frac{2}{n}H_{n-1}B_n(x),
\end{equation}
where $H_{n}=1+\frac{1}{2}+\cdots+\frac{1}{n}$ is the harmonic number and a slight modification of \eqref{4E} gives the identity in \eqref{1A}. 
Let $p(x)=\sum_{k=0}^{n}a_{k}\beta_{k,\lambda}(x)$. Then we have
\begin{align}
a_0&=\bigg\langle \frac{e^{t}-1}{t}\bigg | \frac{\lambda t}{e^{\lambda t}-1}p(x) \bigg \rangle \nonumber\\
&=\frac{2}{n}\sum_{l=0}^{n-2}\frac{1}{n-l}\binom{n}{l}B_{n-l}\bigg\langle \frac{e^{t}-1}{t}\bigg | \frac{\lambda t}{e^{\lambda t}-1}B_{l}(x) \bigg \rangle \label{5E}\\
&\quad\quad+\frac{2}{n}H_{n-1}\bigg\langle \frac{e^{t}-1}{t}\bigg | \frac{\lambda t}{e^{\lambda t}-1}B_{n}(x) \bigg \rangle \nonumber\\
&=\frac{2}{n}\sum_{l=0}^{n-2}\frac{1}{n-l}\binom{n}{l}B_{n-l}\lambda^{l}B_{l}+\frac{2}{n}H_{n-1}\lambda^{n}B_{n}.\nonumber
\end{align}
For $k$, with $ 1 \le k \le n$, we obtain
\begin{align}
k! \lambda^{k-1}a_{k}&=\big \langle \big(e^{\lambda t}-1\big)^{k-1} \big | p(x+1)-p(x)\big \rangle \nonumber\\
&=\frac{2}{n}\sum_{l=1}^{n-2}\frac{1}{n-l}\binom{n}{l}B_{n-l}\big \langle \big(e^{\lambda t}-1\big)^{k-1} \big | B_{l}(x+1)-B_{l}(x)\big \rangle \label{6E}\\
&\quad\quad+\frac{2}{n}H_{n-1}\big \langle \big(e^{\lambda t}-1\big)^{k-1} \big | B_{n}(x+1)-B_{n}(x)\big \rangle \nonumber\\
&=\frac{2}{n}\sum_{l=1}^{n-2}\frac{l}{n-l}\binom{n}{l}B_{n-l}\Delta_{\lambda}^{k-1}0^{l-1}+2 H_{n-1}\Delta_{\lambda}^{k-1}0^{n-1}.\nonumber
\end{align}

Thus, from \eqref{5E} and \eqref{6E} and for $ n \ge 2$, we have
\begin{align*}
&\sum_{k=1}^{n-1}\frac{1}{k(n-k)}B_{k}(x)B_{n-k}(x)\\
&=\frac{2}{n}\bigg(\sum_{l=0}^{n-2}\frac{1}{n-l}\binom{n}{l}B_{n-l}\lambda^{l}B_{l}+H_{n-1}\lambda^{n}B_{n}\bigg)\beta_{0,\lambda}(x)\\
&\quad\quad+\frac{2}{n}\sum_{k=1}^{n}\frac{1}{k!\lambda^{k-1}}\bigg(\sum_{l=1}^{n-2}\frac{l}{n-l}\binom{n}{l}B_{n-l}\Delta_{\lambda}^{k-1}0^{l-1}+n H_{n-1}\Delta_{\lambda}^{k-1}0^{n-1}\bigg)\beta_{k,\lambda}(x).
\end{align*}

(c) In [11], it is shown that the following identity holds for $n \ge 2$:
\begin{equation*}
\sum_{k=1}^{n-1} \frac{1}{k(n-k)}E_{k}(x) E_{n-k} (x) =
\frac{4E_{n+1}}{n^2(n+1)} - \frac{4}{n}\sum_{l=1}^n
\frac{\binom{n}{l}(H_{n-1} - H_{n-l})}{n-l+1} E_{n-l+1} B_{l} (x),
\end{equation*}
where $H_{n}=1+\frac{1}{2}+\cdots+\frac{1}{n}$ is the harmonic number. \par
Write $\sum_{k=1}^{n-1} \frac{1}{k(n-k)}E_{k}(x) E_{n-k} (x)=\sum_{k=0}^{n}a_{k}\beta_{k,\lambda}(x)$. \par
By proceeding similarly to (b), we obtain the following identity:
\begin{align*}
&\sum_{k=1}^{n-1} \frac{1}{k(n-k)}E_{k}(x) E_{n-k} (x) \\
&=\frac{4}{n}\bigg(\frac{E_{n+1}}{n(n+1)} - \sum_{l=1}^n
\frac{\binom{n}{l}(H_{n-1} - H_{n-l})}{n-l+1} E_{n-l+1}\lambda^{l}B_{l}\bigg)\beta_{0,\lambda}(x) \\
&\quad\quad-\frac{4}{n}\sum_{k=1}^{n}\frac{1}{k!\lambda^{k-1}}\bigg(\sum_{l=1}^n
\frac{l\,\binom{n}{l}(H_{n-1} - H_{n-l})}{n-l+1} E_{n-l+1}\Delta_{\lambda}^{k-1}0^{l-1}\bigg)\beta_{k,\lambda}(x).
\end{align*}

(d) In [15], it is proved that the following identity is valid for $n \ge 3$:
\begin{equation}\label{7E}
\sum_{k=1}^{n-1}\frac{1}{k(n-k)}G_k(x)G_{n-k}(x)=-\frac{4}{n}\sum_{k=0}^{n-2}\binom{n}{k}\frac{G_{n-k}}{n-k}B_k(x).
\end{equation}
Again, by proceeding analogously to (b), we get the following identity:
\begin{align*}
\sum_{k=1}^{n-1}\frac{1}{k(n-k)}G_k(x)G_{n-k}(x) &=-\frac{4}{n}\bigg(\sum_{l=0}^{n-2}\binom{n}{l}\frac{G_{n-l}}{n-l}\lambda^{l}B_{l}\bigg)\beta_{0,\lambda}(x)\\
&\quad-\frac{4}{n}\sum_{k=1}^{n-2}\frac{1}{k!\lambda^{k-1}}\bigg(\sum_{l=1}^{n-2}\binom{n}{l}\frac{l\,G_{n-l}}{n-l}\Delta_{\lambda}^{k-1}0^{l-1}\bigg)\beta_{k,\lambda}(x).
\end{align*}

(e) Nielsen [17,2] expressed products of two Bernoulli polynomials in terms of Bernoulli polynomials. Namely, for positive integers $m$ and $n$, with $m+n \ge 2$,
\begin{equation*}
B_{m}(x)B_{n}(x)=\sum_{r}\left\{\binom{m}{2r}n+\binom{n}{2r}m\right\}\frac{B_{2r}B_{m+n-2r}(x)}{m+n-2r}+(-1)^{m+1}\frac{B_{m+n}}{\binom{m+n}{m}}.
\end{equation*}
Again, in a similar way to (b), we can show that
\begin{align*}
&B_{m}(x)B_{n}(x) \\
&=\bigg(\sum_{r}\left\{\binom{m}{2r}n+\binom{n}{2r}m\right\}\frac{B_{2r}B_{m+n-2r}}{m+n-2r}\lambda^{m+n-2r}+(-1)^{m+1}\frac{B_{m+n}}{\binom{m+n}{m}}\bigg)\beta_{0,\lambda}(x) \\
&\quad\quad+\sum_{k=1}^{m+n}\frac{1}{k!\lambda^{k-1}}\sum_{r}\left\{\binom{m}{2r}n+\binom{n}{2r}m\right\}
B_{2r}\Delta_{\lambda}^{k-1}0^{m+n-2r-1}\beta_{k,\lambda}(x).
\end{align*}

(f) Nielsen [17,2] also represented products of two Euler polynomials in terms of Bernoulli polynomials as follows:
\begin{align*}
E_{m}(x)E_{n}(x)&=-2\sum_{r=1}^{m}\binom{m}{r}E_{r} \frac{B_{m+n-r+1}(x)}{m+n-r+1}\\
&\quad-2\sum_{s=1}^{n}\binom{n}{s}E_{s} \frac{B_{m+n-s+1}(x)}{m+n-s+1}\\
&\quad +2(-1)^{n+1} \frac{m!\,n!}{(m+n+1)!} E_{m+n+1}.
\end{align*}
In the same way as (b), we can show that
\begin{align*}
&E_{m}(x)E_{n}(x) \\
&=-2\bigg(\sum_{r=1}^{m}\binom{m}{r}E_{r} \frac{B_{m+n-r+1}\lambda^{m+n-r+1}}{m+n-r+1}+\sum_{s=1}^{n}\binom{n}{s}E_{s} \frac{B_{m+n-s+1}\lambda^{m+n-s+1}}{m+n-s+1} \\
&\quad\quad+(-1)^{n} \frac{m!\,n!}{(m+n+1)!} E_{m+n+1}\bigg)\beta_{0,\lambda}(x) \\
&\quad\quad-\sum_{k=1}^{m+n}\frac{2}{k!\lambda^{k-1}}\bigg(\sum_{r=1}^{m}\binom{m}{r}E_{r}\Delta_{\lambda}^{k-1}0^{m+n-r}+\sum_{s=1}^{n}\binom{n}{s}E_{s}\Delta_{\lambda}^{k-1}0^{m+n-s}\bigg)\beta_{k,\lambda}(x).
\end{align*}

(g) As the last example, we would like to express the sum on the left hand side of \eqref{7E} in terms of the degenerate Bernoulli polynomials of order $r$. For this, we first observe that, for any positive integer $a$, 
\begin{equation}\label{8E}
\bigg(\frac{e^{t}-1}{t}\bigg)^{a}B_{n}^{(r)}\big(\frac{x}{\lambda}\big)=I^{a}B_{n}^{(r)}\big(\frac{x}{\lambda}\big)=\frac{\lambda^{a}}{\langle n+1 \rangle_{a}}\sum_{m=0}^{a}(-1)^{a-m}\binom{a}{m}
B_{n+a}^{(r)}\big(\frac{x+m}{\lambda}\big),
\end{equation}
where $I$ is the integral operator $Iq(x)=\int_{x}^{x+1}q(u)\,du,$ and $\langle n+1 \rangle_{a}=(n+1)\cdots(n+a)$.
Let $p(x)=\sum_{k=1}^{n-1}\frac{1}{k(n-k)}G_k(x)G_{n-k}(x),\,\,(n \ge 3)$. Then $p(x)=-\frac{4}{n}\sum_{l=0}^{n-2}\binom{n}{l}\frac{G_{n-l}}{n-l}B_l(x).$
In order to apply Theorem 4.2, we need to compute (see \eqref{13D}) \par
\begin{equation*}
I^{r-k} p\big(\lambda {\bf{B}}^{(r-k)}\big(\frac{x}{\lambda}\big)\big)  \bigg |_{x=j}
= \bigg \langle e^{jt} \bigg(\frac{e^{t}-1}{t}\bigg)^{r-k} \bigg| \bigg(\frac{\lambda t}{e^{\lambda t}-1} \bigg)^{r-k} p(x) \bigg \rangle,
\end{equation*}
which is equal to
\begin{align}
&I^{r-k} p\big(\lambda {\bf{B}}^{(r-k)}\big(\frac{x}{\lambda}\big)\big)  \bigg |_{x=j} \nonumber\\
&=-\frac{4}{n}\sum_{l=0}^{n-2}\binom{n}{l}\frac{G_{n-l}}{n-l}\bigg \langle e^{jt} \bigg(\frac{e^{t}-1}{t}\bigg)^{r-k} \bigg| \bigg(\frac{\lambda t}{e^{\lambda t}-1} \bigg)^{r-k} B_l(x) \bigg \rangle \label{9E}\\
&=-\frac{4}{n}\sum_{l=0}^{n-2}\binom{n}{l}\frac{G_{n-l}}{n-l}\lambda^{l} \bigg \langle e^{jt}  \bigg| \bigg(\frac{e^{t}-1}{t}\bigg)^{r-k-1}B_{l}^{(r-k)}\big(\frac{x}{\lambda}\big) \bigg \rangle \nonumber \\
&=-\frac{4}{n}\sum_{l=0}^{n-2}\binom{n}{l}\frac{G_{n-l}}{n-l}\frac{\lambda^{l+r-k-1}}{\langle l+1 \rangle_{r-k-1}}\sum_{m=0}^{r-k-1}(-1)^{r-k-1-m}\binom{r-k-1}{m}B_{l+r-k-1}^{(r-k)}\big(\frac{j+m}{\lambda}\big), \nonumber
\end{align}
where we used \eqref{8E}.
Assume that $n \ge r$ and $ n \ge 3$. Then, from Theorem 4.2 and \eqref{9E}, we obtain
\begin{align*}
\sum_{k=1}^{n-1}\frac{1}{k(n-k)}&G_k(x)G_{n-k}(x)=
-\frac{4}{n}\sum_{k=0}^{r-1}\bigg(\frac{1}{k!}\sum_{j=0}^{k}\sum_{l=0}^{n-2}\sum_{m=0}^{r-k-1}(-1)^{r-j-m-1}\binom{k}{j}\binom{n}{l}\\
&\times \frac{G_{n-l}}{n-l}\frac{\lambda^{l+r-k-1}}{\langle l+1 \rangle_{r-k-1}}\binom{r-k-1}{m}B_{l+r-k-1}^{(r-k)}\big(\frac{j+m}{\lambda}\big)\bigg) \beta_{k,\lambda}^{(r)}(x) \\
&+\sum_{k=r}^{n-2}\bigg(\frac{1}{k!\lambda^{k-r}}\sum_{j=0}^{r}\sum_{l=0}^{k-r}\sum_{m=1}^{n-1}\frac{(-1)^{k-j-l}}{{m(n-m)}}\binom{r}{j}\binom{k-r}{l}G_{m}(j+l \lambda)\\
&\times G_{n-m}(j+l \lambda)\bigg)\beta_{k,\lambda}^{(r)}(x). 
\end{align*}
\section{Conclusion}

In this paper, we were interested in representing any polynomial in terms of the degenerate Bernoulli polynomials and of the higher-order degenerate Bernoulli polynomials. We were able to derive formulas for such representations with the help of umbral calculus. We showed that,  by letting $\lambda$ tends to zero, they agree with the previously found formulas for representations by the Bernoulli polynomials and by the higher-order Bernoulli polynomials. Further, we illustrated the formulas with some examples. \par
Even though the method adopted in this paper is elementary, they are very useful and powerful. Indeed, as we mentioned in the Section 1, both a variant of Miki's identity and Faber-Pandharipande-Zagier (FPZ) identity follow from the one identity (see \eqref{1A}) that can be derived from a formula (see Remark 3.2) involving only derivatives and integrals of the given polynomial, while all the other proofs are quite involved. We recall here that the FPZ identity was a conjectural relations between Hodge integrals in Gromov-Witten theory. In addition, this paper demonstrates that umbral calculus is a very convenient tool.\par
It is one of our future research projects to continue to find formulas representing polynomials in terms of some specific special polynomials and to apply those in discovering some interesting identities.


\end{document}